\documentclass[12pt]{article}
\usepackage{amsmath,amsfonts,amssymb,theorem}
\textheight215truemm \textwidth 138truemm  \date{}

 \newcommand{\C}{\mathbb C}
 \newcommand{\PP}{\mathbb P}
 \newcommand{\Z}{\mathbb Z}
 \newcommand{\hX}{\widehat X}
 \newcommand{\SL}{\operatorname{SL}}

\begin{document}
 \title{Warwick EuroConference (July 1996)}
 \author{Klaus Hulek and others (Editors)}
 \date{}
 \maketitle

 \begin{abstract}
 This file, provided as bibliographical information for the benefit of
contributors, contains the table of contents and editors' foreword from the
proceedings of the July 1996 Warwick EuroConference on Algebraic Geometry.
 \end{abstract}

 \subsection*{Full title:} Recent Trends in Algebraic Geometry --
EuroConference on Algebraic Geo\-metry (Warwick, July 1996), Editors: Klaus
Hulek (chief editor), Fabrizio Catanese, Chris Peters and Miles Reid, CUP,
May 1999

 \subsection*{Contents:} Victor V. Batyrev: Birational Calabi--Yau $n$-folds
have equal Betti numbers, 1--11 (alg-geom/9710020)

{\parindent0pt

Arnaud Beauville: A Calabi--Yau threefold with non-Abelian fundamental group,
13--17 (alg-geom/9502003)

K. Behrend: Algebraic Gromov--Witten invariants, 19--70 (compare
alg-geom/9601011)

Philippe Eyssidieux: K\"ahler hyperbolicity and variations of Hodge structures,
71--92

Carel Faber: Algorithms for computing intersection numbers on moduli spaces of
curves, with an application to the class of the locus of Jacobians, 93--109
(alg-geom/9706006)

Marat Gizatullin: On some tensor representations of the Cremona group of the
projective plane, 111--150

Y. Ito and I. Nakamura: Hilbert schemes and simple singularities, 151--233

Oliver K\"uchle and Andreas Steffens: Bounds for Seshadri constants, 235--254
(alg-geom/9601018)

Marco Manetti: Degenerate double covers of the projective plane, 255--181 (see
further math.AG/9802088)

David R. Morrison: The geometry underlying mirror symmetry, 283--310
(alg-geom/9608006)

Shigeru Mukai: Duality of polarized K3 surfaces, 311--326

Roberto Paoletti: On symplectic invariants of algebraic varieties coming from
crepant contractions, 327--346

Kapil H. Paranjape: The Bogomolov--Pantev resolution, an expository account,
347--358 (math.AG/9806084)

Tetsuji Shioda: Mordell--Weil lattices for higher genus fibration over a
curve, 359--373

Bernd Siebert: Symplectic Gromov--Witten invariants, 375--424

Claire Voisin: A generic Torelli theorem for the quintic threefold, 425--463

P.M.H. Wilson: Flops, Type~III contractions and Gromov--Witten invariants on
Calabi--Yau threefolds, 465--483

}

 \section*{Foreword}
 \markboth{\qquad Foreword\hfill}{\hfill Foreword\qquad}

The volume contains a selection of seventeen survey and research articles
from the July 1996 Warwick European algebraic geometry conference. These
papers give a lively picture of current research trends in algebraic
geo\-metry, and between them cover many of the outstanding hot topics in
the modern subject. Several of the papers are expository accounts of
substantial new areas of advance in mathematics, carefully written to be
accessible to the general reader. The book will be of interest to a wide
range of students and nonexperts in different areas of mathematics,
geometry and physics, and is required reading for all specialists in
algebraic geometry.

The European algebraic geometry conference was one of the climactic events
of the 1995--96 \hbox{EPSRC} Warwick algebraic geometry symposium, and
turned out to be one of the major algebraic geometry events of the 1990s.
The scientific committee consisted of A.~Beauville (Paris), F.~Catanese
(Pisa), K.~Hulek (Hannover) and C.~Peters (Grenoble) representing AGE
(\hbox{Algebraic} Geometry in Europe, an EU HCM--TMR network) and
N.J.~Hitchin (Oxford), J.D.S.~Jones and M.~Reid (Warwick) representing
Warwick and British mathematics. The conference attracted 178 participants
from 22 countries and featured 33 lectures from a star-studded cast of
speakers, including most of the authors represented in this volume.

\paragraph{The expository papers}
 Five of the articles are expository in intention: among these a beautiful
short exposition by Paranjape of the new and very simple approach to the
resolution of singularities; a detailed essay by Ito and Nakamura on the
ubiquitous ADE classification, centred around simple surface
singularities; a discussion by Morrison of the new special Lagrangian
approach giving geometric foundations to mirror symmetry; and two deep and
informative survey articles by Behrend and Siebert on Gromov--Witten
invariants, treated from the contrasting viewpoints of algebraic and
symplectic geometry.

\paragraph{Some main overall topics}
 Many of the papers in this volume group around a small number of main
research topics. Gromov--Witten invariants was one of the main new
breakthroughs in geometry in the 1990s; they can be developed from several
different starting points in symplectic or algebraic geometry. The survey
of Siebert covers the analytic background to the symplectic point of view,
and outlines the proof that the two approaches define the same invariants.
Behrend's paper explains the approach in algebraic geo\-metry to the
invariants via moduli stacks and the virtual fundamental class, which
essentially amounts to a very sophisticated way of doing intersection
theory calculations. The papers by Paoletti and Wilson give parallel
applications of Gromov--Witten invariants to higher dimensional varieties:
Wilson's paper determines the Gromov--Witten invariants that arise from
extremal rays of the Mori cone of Calabi--Yau 3-folds, whereas Paoletti
proves that Mori extremal rays have nonzero associated Gromov--Witten
invariants in many higher dimensional cases. The upshot is that extremal
rays arising in algebraic geometry are in fact in many cases invariant in
the wider symplectic and topological setting.

 Another area of recent spectacular progress in geometry and theoretical
physics is Calabi--Yau 3-folds and mirror symmetry. This was another major
theme of the EuroConference that is well represented in this volume. The
paper by Voisin, which is an extraordinary computational tour-de-force,
proves the generic Torelli theorem for the most classical of all
Calabi--Yau 3-folds, the quintic hypersurface in $\PP^4$. The survey by
Morrison explains, among other things, the Strominger--Yau--Zaslow special
Lagrangian interpretation of mirror symmetry. Beauville's paper gives the
first known construction of a Calabi--Yau 3-fold having the quaternion
group of order 8 as its fundamental group. The paper by Batyrev proves
that the Betti numbers of a Calabi--Yau 3-fold are birationally invariant,
using the methods of $p$-adic integration and the Weil conjectures; the
idea of the paper is quite startling at first sight (and not much less so
at second sight), but it is an early precursor of Kontsevich's idea of
motivic integration, as worked out in papers of Denef and Loeser. Several
other papers in this volume (those of Ito and Nakamura, Mukai,
\hbox{Shioda} and Wilson) are implicitly or explicitly related to
Calabi--Yau 3-folds in one way or another.

\paragraph{Other topics}
 The remaining papers, while not necessarily strictly related in subject
matter, include some remarkable achievements that illustrate the breadth
and depth of current research in algebraic geometry. Shioda extends his
well-known results on the Mordell--Weil lattices of elliptic surfaces to
higher genus fibrations, in a paper that will undoubtedly have substantial
repercussions in areas as diverse as number theory, classification of
surfaces, lattice theory and singularity theory. Faber continues his study
of tauto\-logical classes on the moduli space of curves and Abelian
varieties, and gives an algorithmic treatment of their intersection
numbers, that parallels in many respects the Schubert calculus; he obtains
the best currently known partial results determining the class of the
Schottky locus. Gizatullin initiates a fascinating study of
representations of the Cremona group of the plane by birational
transformations of spaces of plane curves. Eyssidieux gives a study, in
terms of Gromov's K\"ahler hyperbolicity, of universal inequalities
holding between the Chern classes of vector bundles over Hermitian
symmetric spaces of noncompact type admitting a variation of Hodge
structures. K\"uchle and Steffens' paper contains new twists on the idea
of Seshadri constants, a notion of local ampleness arising in recent
attempts on the Fujita conjecture; they use in particular an ingenious
scaling trick to provide improved criteria for the very ampleness of
adjoint line bundles.

Manetti's paper continues his long-term study of surfaces of general type
constructed as iterated double covers of $\PP^2$. He obtains many
constructions of families of surfaces, and proves that these give complete
connected components of their moduli spaces, provided that certain
naturally occuring degenerations of the double covers are included. This
idea is used here to establish a bigger-than-polynomial lower bound on the
growth of the number of connected components of moduli spaces. In more
recent work, he has extended these ideas in a spectacular way to exhibit
the first examples of algebraic surfaces that are proved to be
diffeomorphic but not deformation equivalent.

The Fourier--Mukai transform is now firmly established as one of the most
important new devices in algebraic geometry. The idea, roughly speaking,
is that a sufficiently good moduli family of vector bundles (say) on a
variety $X$ induces a correspondence between $X$ and the moduli space
$\hX$. In favourable cases, this correspondence gives an equivalence of
categories between coherent sheaves on $X$ and on $\hX$ (more precisely,
between their derived categories). The model for this theory is provided
by the case originally treated by Mukai, when $X$ is an Abelian variety
and $\hX$ its dual; Mukai named the transform by analogy with the
classical Fourier transform, which takes functions on a real vector space
to functions on its dual. It is believed that, in addition to its many
fruitful applications in algebraic geometry proper, this correspondence
and its generalisations to other categories of geometry will eventually
provide the language for mathematical interpretations of the various
``dualities'' invented by the physicists, for example, between special
Lagrangian geometry on a Calabi--Yau 3-fold and coherent algebraic
geometry on its mirror partner (which, as described in Morrison's article,
is conjecturally a fine moduli space for special Lagrangian tori). Mukai's
magic paper in this volume presents a Fourier--Mukai transform for K3
surfaces, in terms of moduli of semi-rigid sheaves; under some minor
numerical assumptions, he establishes the existence of a dual K3 surface,
the fact that the Fourier--Mukai transform is an equivalence of derived
categories, and the biduality result in appropriate cases.

The paper of Ito and Nakamura is the longest in the volume; it combines a
detailed and wide-ranging expository essay on the ADE classification with
an algebraic treatment of the McKay correspondence for the Kleinian
quotient singularities $\C^2/G$ in terms of the $G$-orbit Hilbert scheme.
The contents of their expository section will probably come as a surprise
to algebraic geo\-meters, since alongside traditional aspects of simple
singularities and their ADE homologues in algebraic groups and
representation of quivers, they lay particular emphasis on partition
functions in conformal field theory with modular invariance under
$\SL(2,\Z)$ and on II$_1$ factors in von Neumann algebras. Their study of
the $G$-Hilbert scheme makes explicit for the first time many aspects of
the McKay correspondence relating the exceptional locus of the Kleinian
quotient singularities $\C^2/G$ with the irreducible representations of
$G$; for example, the way in which the points of the minimal resolution
can be viewed as defined by polynomial equations in the character spaces
of the corresponding irreducible representations, or the significance in
algebraic terms of tensoring with the given representation of $G$. Ito and
Nakamura and their coworkers are currently involved in generalising many
aspects of the $G$-orbit Hilbert scheme approach to the resolution of
Gorenstein quotient singularities and the McKay correspondence to finite
subgroups of $\SL(3,\C)$, and this paper serves as a model for what one
hopes to achieve.

\paragraph{Thanks to all our sponsors} 
The principal financial support for the Euro\-Conference was a grant of 
ECU40,000 from EU TMR (Transfer and Mobility of Researchers), contract
number ERBFMMACT 950029; we are very grateful for this support, without
which the conference could not have taken place. The main funding for the
1995--96 Warwick algebraic geometry symposium was provided by British
EPSRC (Engineering and physical sciences research council). Naturally
enough, the symposium was one of the principal activities of the Warwick
group of AGE (European Union HCM project Algebraic Geometry in Europe,
Contract number ERBCHRXCT 940557), and financial support from Warwick AGE
and the other groups of AGE was a crucial element in the success of the
symposium and the EuroConference. We also benefitted from two visiting
fellowships for Nakamura and Klyachko from the Royal Society (the UK
Academy of Science). Many other participants were covered by their own
research grants.

The University of Warwick, and the Warwick Mathematics Institute also
provided substantial financial backing. All aspects of the conference were
enhanced by the expert logistic and organisational help provided by the
\hbox{Warwick} Math Research Centre's incomparable staff, Elaine Greaves
Coelho, Peta McAllister and Hazel Graley.

\bigskip

 \centerline{Klaus Hulek and Miles Reid, November 1998}


 \end{document}